\documentclass[a4paper,12pt]{article}
\usepackage[latin1]{inputenc}
\usepackage{amsmath,amssymb,amsthm,cite,verbatim}
\usepackage[all]{xy}
\usepackage[bookmarks=true]{hyperref}
\usepackage{color}
\usepackage{graphicx}

\title{Hyperelliptic parametrizations of $\Q$-curves}
\author{ \ Francesc Bars\footnote{First and third author are supported by MTM2016-75980-P and MDM-2014-0445.}, Josep Gonz\'alez  \footnote{Second author is
partially supported by DGI grant  MTM2015-66180-R.} and Xavier Xarles}

\newtheorem{prop}{Proposition}
\newtheorem{cor}{Corollary}
\newtheorem{rem}{Remark}

\theoremstyle{definition}

\theoremstyle{remark}

\newtheorem*{exemp}{Example}

\numberwithin{equation}{section}

\newcommand{\Q}{\mathbb{Q}}

\newcommand{\Z}{\mathbb{Z}}

\newcommand{\bP}{\mathbb{P}}

\newcommand{\C}{\mathbb{C}}

\newcommand{\Gal}{\mathrm{Gal}}

\newcommand{\Aut}{\operatorname{Aut}}

\newcommand{\Jac}{\operatorname{Jac}}

\newcommand{\dv}{\operatorname{div}}

\newcommand{\rank}{\operatorname{rank}}

\newcommand{\cL}{{\mathcal L}}

\date{}
\begin{document}
\maketitle

\begin{abstract}
For a square-free integer $N$, we present a procedure to compute $\Q$-curves parametrized by rational points of the modular curve $X_0^*(N)$ when this is hyperelliptic.
\end{abstract}

\section{Introduction}
Let $X$ be a curve defined over $\overline{\Q}$. The curve $X$ is said to be a $\Q$-curve if it is isogenous to all Galois conjugates.  In \cite{Elkies}, Elkies proved that every $\Q$-curve without complex multilication (CM) is isogenous over $\overline{\Q}$ to a $\Q$-curve attached to a rational point of the modular curve $X_0^*(N)=X_0(N)/B(N)$ for some square-free integer $N$, where $B(N)$ denotes the group of the Atkin-Lehner involutions of $X_0(N)$. Every rational non-cuspidal point in $X_0^*(N)$ lifts to $X_0(N)$ providing $\Q$-curves, with or without CM, defined over abelian extensions of $\Q$ of type $(2,\cdots ,2)$.

In \cite{GL}, it is given a procedure to parametrize the $j$-invariants of these $\Q$-curves when the genus $g_N^*$ of $X_0^*(N)$ is at most $1$. In these cases, the set $X_0^*(N)(\Q)$ is infinite. Basically, in this paper it is given a method to determine the symmetric functions of the set  $\{ j(d\,z)\colon 1\leq d|N\}$, where $j(z)$ denotes the usual generator of $\Q(X_0(1))$,  from a suitable  generators of $\Q(X_0^*(N))$. Here, we present a similar procedure for the case that $X_0^*(N)$ is hiperelliptic, that amounts to saying $g_N^*=2$. In fact, there are exactly $36$ values of $N$ for which $g_N^*=2$. By Faltings, for such a value of $N$, we are dealing with a finite number of cases.
 Firstly, we consider the rational points provided by Magma.
 Although the rank of $\Jac (X_0^*(N))$ is equal to $2$ and the classical  Chabauty method does not work, we can  determine the full set $X_0^*(N)(\Q)$   for $19$ values of $N$ by using a Chabauty procedure on a finite set of unramified $2$-coverings of the curve.

The article is organized as follows.  In Proposition \ref{main} of \S 2, we present the main tool to parametrize $\Q$-curves from rational  non-cuspidal  points of $X_0^*(N)$. In \S 3,
we give a list of equations of $X_0^*(N)$ when $g_N^*=2$ together their rational points provided by Magma  and, in Proposition \ref{rational}, we determine all rational points  for $19$ values of $N$.
 In \S 4, we show how the $j$-invariants of
 the $\Q$-curves  curves over these rational points are computed for the case  $N=67$. Next, for  all values of $N$  we determine which of the parametrized $\Q$-curves  have CM and, for all of them, we give the discriminant  $D$ of the  order of its endomorphism ring.     Moreover, if the $j$-invariant of the $\Q$-curve lies in a quadratic field, it is given explicitly;     otherwise,
 we provide the   number field  $\Q(j       )$.

 We recall     that there is a finite number of discriminants $D$ of orders of imaginary quadratic
 fields $K$ such that     $\Gal(H_D/\Q)$      is of the type      $(2,\cdots ,2)$, where $H_D$ denotes the ring class field of the quadratic order of
  discriminant $D$. In fact, this condition is equivalent to say that the $j$-invariant of an  elliptic curve with CM by the order of discriminant $D$  generates a totally real number field. Moreover,  $|\Gal(H_D/K)| $ divides $16$ (for more detail, see \cite{Buell} and CM-computations in the web-page https://mat-web.upc.edu/people/joan.carles.lario/).
  The results obtained  when $g_N^*=2$ show that
  almost these  $\Q$-curves have CM, which reinforces the conjecture that for a large enough $N$, the curve $X_0^*(N)$ does not have rational points
  parametrizing $\Q$-curves without CM.
\
\section{Preliminary results}
Let $X$ be a genus two curve defined over a subfield $K$ of the complex field $\C$ that is the normalization of the curve given by the affine equation
$$
y^2=x^6+a_5 x^5+a_4x^4+a_3x^3+a_2x^2+a_1x+a_0\,,
$$
with $a_i\in K$ for all $i \leq 5$. Let us denote by $w$ the hyperelliptic involution. There are two points $P_0$ and $P_1=w(P_0)$ over the singularity at infinity, defined over $K$ that  are not Weierstrass points and both  satisfy $\displaystyle{\left(\frac{y}{x^3}\right)^2(P_i)=1}$. Denote by  $P_0$ such a point satisfying $\displaystyle{\left(\frac{y}{x^3}\right)(P_0)=1}$ and, thus, $\displaystyle{\left(\frac{y}{x^3}\right)(P_1)=-1}$.
For an integer $n\geq 1$, we consider the  $K$-vector space
$$
\cL_i=\{ f\in K(X)\colon \dv (f)\geq -n \,P_0\}\,.
$$
It is clear that   $\dim \cL_n=1$ for $n\leq 2$ and, by the Riemann-Roch Theorem, we know that  $\dim \cL_n=n-1$ when $n>2$. We denote by $\dv^- f$ the polar part of $\dv f$, i.e. $\dv f+.\dv^- f\geq 0$.
\begin{prop}\label{main}
Keeping the above notation, the functions $f_3,f_4,f_5\in K(X)$ defined by
$$
\begin{array}{cr}
f_3=&\displaystyle{\frac{8 a_3 - 4 a_4\,a_5+ a_5^3}{32} +\frac{4 a_4 - a_5^2}{16} x + \frac{a_5}{4} x^2 +\frac{1}{2} x^3 +\frac{1}{2}y\,,}\\
f_4=& \displaystyle{
\frac{64 a_2 - 16 a_4^2 - 32 a_3 a_5 + 24 a_4 a_5^2 - 5 a_5^4}{256}+
x\,f_3\,,}\\
f_5=&\displaystyle{
\frac{128 a_1 - 64 a_3 a_4 - 64 a_2 a_5 + 48 a_4^2 a_5 + 48 a_3 a_5^2 -   40 a_4 a_5^3 + 7 a_5^5}{512}+ x\,f_4\,,
    }  \\
\end{array}
$$
 vanish at $w(P_0)$ and satisfy that  $\dv^{-}f_i=i\,P_0$ for all $i\in\{3,4,5\}$. In particular,  $\cL_3=\langle 1,f_3\rangle$, $\cL_4=\langle 1,f_3,f_4\rangle$  and $\cL_5=\langle 1,f_3,f_4,f_5\rangle$.

\end{prop}
\noindent{\bf Proof.} Let $g_3\in \cL_3$ be a non constant function. By adding a constant, if necessary, we can assume $g_3(w(P_0))=0$. We have that the function $h=g_3-(g_3|w)$ satisfies $h|w=-h$ and $\dv ^-h=3(P_0)+3 (w(P_0))$. Hence, $h=A\, y$ for some non zero $A\in K$. Putting $f_3= g_3/A$, we get  $f_3-(f_3|w)=y$. Moreover, the functions $f_3 + (f_3|w)$ and $f_3\cdot (f_3|w)$ are invariant under the action of $w$ and satisfy
 $$\dv^-(f_3 + (f_3|w))= 3(P_0)+3 (w(P_0))\quad \text{ and }\quad \dv^- (f_3\cdot (f_3|w))\leq 2(P_0)+2 (w(P_0))\,.$$ Therefore, $f_3 + (f_3|w)$ is a polynomial in $x$ of degree $3$ and
 $f_3\cdot (f_3|w) $ must be a polynomial in $x$ of degree at most $2$. Since
$$x^6+a_5 x^5+a_4x^4+a_3x^3+a_2x^2+a_1x+a_0-(f_3+(f_3|w))^2=-4((f_3\cdot (f_3|w))\,,$$ the polynomial $P(X)^2=(f_3+(f_3|w))^2$ is determined and, thus, $P(x)$ is determined  up to sign. Hence, $f_3$ must be $1/2( y \pm P(x))$. Since $(y/x^3)(P_0)=1$, we take the sign such that $f_3=1/2(y+ x^3+\cdots)$.

 The function $x f_3$ lies in $\cL_4$ and   $\dv ^- x\,f_3=4 P_0$. There exists $k\in K$ such that the function  $f_4=x\,f_3+k$ vanishes at $w(P_0)$. By construction,   $f_4-(f_4|w)=x\cdot y$ and  $f_4+(f_4|w)$ is a polynomial in $x$ of degree $4$.
We can determine $k$ by using  that the function  $x^2 y^2 -(f_4+(f_4|w))^2$ is a polynomial in $x$ of degree at most $3$, because it is equal to the function   $ -4((f_4\cdot (f_4|w))$.

Similarly, $x\,f_4\in\cL_5\backslash \cL_4$ and there exists $k\in K$ such that the function   $f_5=x\,f_4+k$ vanishes at $w(P_0)$. Now,  $f_5-(f_5|w)=x^2\cdot y$ and  $f_5+(f_5|w)$ is a polynomial in $x$ of degree $5$. We  determine $k$ by using
that  $x^4 y^2 -(f_5+(f_5|w))^2$  is a polynomial in $x$ of degree at most $4$.
\hfill $\Box$
\begin{cor}
For an integer $n\geq 3$ and  function $f\in \cL_n$, the function $f-f(\infty')$ is a $K$-linear combination of the  functions $\{ f_5f_3^k,f_4 f_3^k, f_3^{k+1}\colon 0\leq  k\leq \lfloor n/3\rfloor\}$.
\end{cor}
\noindent{\bf Proof.} For an integer $n\geq 3$, let $i\in\{0,1,2\}$ be such that $n\equiv i\pmod 3$. The statement follows from the fact that  the function $h=f_{i+3}f_3^{ (n-i)/3-1}$ lies in $\cL_n$ with  $\dv ^-h =n P_0$ and $h(\infty')=0$.
\hfill $\Box$
\begin{rem} Note that if $q$ is an analytic  uniformizing parameter at $P$ such that $x=1/q+\cdots$ and $y=1/q^3+\cdots$, then $f_i=1/q^i+\cdots$ for $ i \in \{3,4,5\}$.
\end{rem}
\section{Application to genus two curves $X_0^*(N)$} In \cite{HaHa}, it is proved that when $N$ is square-free,  $X_0^*(N)$ is hyperelliptic if, and only if, it  has genus two. There are $35$ square-free integers $N$ such that $X_0^*(N)$ has genus two (cf. \cite[Remark 1]{Ha97}). In all these cases, there is an only basis $h_1$ and $h_2$ of $S_2(\Gamma_0(N))^{B(N)}$ such that  their  $q$-expansions lie in $\Z[[q]]$ and are of the form $h_1(q)=q+\sum_{n\geq 3} b_n q^n$ and $h_2(q)=q^2+\sum_{n\geq 3}c_n q^n$.
The functions on $X_0^*(N)$ defined as follows
$$
x= \frac{h_1}{h_2}=1/q+\cdots \,,\quad y=- q\frac{ d\, x}{d\,q}/h_2=1/q^3+\cdots
$$
satisfy an equation of the form $
y^2=x^6+a_5 x^5+a_4x^4+a_3x^3+a_2x^2+a_1x+a_0$ with $a_i\in\Z$ for all $i$. Denoting by $\infty$ the infinity cusp and by $\infty'=w(\infty)$, we have $\dv^{-}x= \infty+\infty'$,  $\dv^- y=3\infty+3 \infty'$ and $(y/x^3)(\infty)=1$.

\vskip 0.2 cm
Consider the symmetric functions obtained from the functions $j(d \, z)$ for $1\leq d|N$:
$$J_1(z)=\sum_{1\leq d|N}j(d\,z)\,,\cdots\,,J_m(z)=\prod_{1\leq d|N}j(d\,z)\,,$$
where $m=2^{\omega(N)}$ and  $\omega (N)$ denotes the number of primes dividing $N$.
We determine every function $J_i$ as a  $\Q$-linear combination of functions of the form $ f_5f_3^k,f_4 f_3^k, f_3^{k}$ for $ k\geq    0$. Given a non-cuspidal $Q\in X_0^*(N)$, the $j$-invariants of the   $\Q$-curves attached to this point  are the roots of the polynomial in $z$:
$$
z^m + \sum_{i=1}^{n}(-1)^iJ_{i}( Q) z^{n-i}\,.
$$

  We know that for a non trivial automorphism $u$ of $X_0^*(N)$, one has  $u(\infty)\neq \infty$ (cf. \cite[Lemma 3.1]{BH}). Hence,  for all these curves $|X_0^*(N)(\Q)\backslash\{\infty\}|\geq 1$ and for the bielliptic curves, i.e. for $N\in \{106,122,129, 158,166,215,390\}$
  (cf. \cite[Theorem 1]{BaGon}),  we have that $|\Aut (X_0^*(N))|=4$ and, thus,   $|X_0^*(N)(\Q)\backslash\{\infty\}|\geq 3$.
  Next, in Table 1, we present the equations  with the functions $x$ and $y$  obtained  following the procedure  mentioned above,
  together with the rational points with $x$-coordinate with height less than $10^4$.

  \newpage
\subsection{Equations and rational points}
\begin{small}
$$
\begin{array}{c|l|l|}
 N& \text{equation}& X_0^*(N)(\Q)\backslash\{\infty,\infty'\}\\\hline
67 &y^2=x^6 - 4 x^5 + 6 x^4 - 6 x^3 + 9 x^2 - 14 x + 9 & (-1,\pm 7), (0,\pm 3),(1,\pm 1), (2,\pm 1)\\
73 &y^2=   x^6- 4 x^5+ 6 x^4 + 2 x^3 - 15 x^2  + 10 x +1 & (0,\pm 1), (1 ,\pm 1),(2,\pm 3), (\frac{3}{2},\pm \frac{5}{8})\\
85 &y^2=(x^2-2x+5) (    x^4- 2x^3+ 3 x^2- 6 x+5) &(0,\pm 5),(1,\pm 2), (2,\pm 5), (\frac{3}{2},\pm \frac{17}{8}), (-\frac{4}{3},\pm \frac{425}{27}) \\
93 & y^2=( x^3-2 x^2-x+3) ( x^3+2 x^2 -5x+3 )& (-1,\pm 3), (0,\pm 3),(1,\pm 1), (2,\pm 3), (\frac{3}{2}, \pm \frac{9}{8}),\\ & &(\frac{1}{4},\pm \frac{143}{64})\\ 103 & y^2=x^6 - 10 x^4+ 22 x^3- 19 x^2+ 6 x  +1  &  (0,\pm 1), (1,\pm 1), (3,\pm 19)\\
106 &y^2= x^6-4x^5+4x^4+2x^3+4x^2-4x+1 & (-1,\pm 4), (0,\pm 1), (1,\pm 2), (2,\pm 5),(\frac{1}{2},\pm \frac{5}{8}) \\
107 & y^2=x^6 -4 x^5 +10 x^4 -18 x^3+17 x^2 -10 x+1& (0,\pm 1),(2, \pm 1) \\
 115 &y^2=( x^3-2x^2+3x-1) ( x^3+2 x^2-9x+7) & (1,\pm 1),(2,\pm 5),(\frac{1}{2},\pm \frac{5}{8}), (\frac{4}{3},\pm \frac{35}{27})\\
 122 &y^2=  x^6+4x^4-6x^3+4x^2+1& (-1, \pm 4), (0, \pm 1), (1, \pm 2), (\frac{3}{2}, \pm \frac{37}{8}),( \frac{2}{3}, \pm \frac{37}{27} )\\
 129 &y^2= x^6-4x^5-4x^4+12 x^3+4x^2-12x+4 & (-1, \pm 3), (0, \pm 2), (1, \pm 1), (\frac{1}{2},\pm \frac{3}{8}), \\ & & (-\frac{7}{ 5},\pm  \frac{383}{125}), (\frac{7}{12}, \pm \frac{383}{1728})\\
 133 &y^2= x^6+ 4x^5-18 x^4+26 x^3-15 x^2+2x+1 & (0, \pm 1), (1,\pm 1), (\frac{3}{5},\pm \frac{ 83}{125})\\
 134 &y^2=x^6-4 x^5+2 x^4-2 x^3+x^2+2x+1  & (-1, \pm 3), (0, \pm 1), (1, \pm 1), (-\frac{1}{2}, \frac{7}{8})\\
 146 &y^2= x^6 -4 x^5+2x^4+6x^3+x^2+2x+1 & (-1, \pm 1), (0, \pm 1), (1, \pm 3), (2, \pm 5)\\
 154 &y^2= (x-2) (x^2+x+2) ( x^3-3 x^2-x-1)& (0, \pm 2), (1, \pm 4), (2, 0), (-\frac{3}{2},\pm \frac {77}{8}),  \\ & &(-\frac{1}{3}, \pm\frac{56}{27}), (4, \pm  22)\\
 158 &y^2=x^6-4x^4+2x^3-4x^2+1&(0, \pm 1), (2,\pm  1), (\frac{1}{2}, \pm \frac {1}{8})\\
161 & y^2=( x^3-2 x^2+3x-1)(x^3+2x^2+3x-5)& (-1,\pm  7), (1, \pm 1), (-\frac{1}{2}, \pm \frac{35}{8}), (-\frac{1}{4},\pm \frac{ 209}{64})\\
165 &y^2=
(x-1)(x+3)(x^2-x-1)( x^2-x+3) & (-1, \pm 4), (0, \pm 3),(1,0), (2,\pm 5), (-\frac{1 }{2},\pm \frac{15}{8}),\\ &  &(-3,0), (\frac{2}{3},\pm \frac{55}{27}), (\frac{5}{2}, \pm \frac{99}{8})\\
167 & y^2=x^6-4 x^5+2 x^4-2 x^3-3 x^2+2x-3 & (-1,\pm 1)\\
170 & y^2=( x^2-5 x+5) (x^4-11 x^3+48 x^2 -87x +53)& (1,  \pm 2), (2, \pm 1), (\frac{3}{2}, \pm \frac{5}{8}), (4, \pm  5), (\frac{11}{3}, \frac{38}{27})\\
177 &y^2= x^6+2 x^4-6 x^3+5x^2-6x+1& (0,\pm 1), (\frac{3}{2}, \pm \frac{17}{8}) \\
 186 &y^2= ( x^3-2x^2+ x+1) ( x^3+2 x^2+5x+1) & (-1, \pm 3), (0, \pm 1), (1, \pm 3),  (2,\pm 9),(-\frac{1}{2},\pm \frac {3}{8}),\\  & & (-\frac{4}{3},\pm  \frac{143}{27})\\
 191 & y^2= x^6+2 x^4+2x^3 +5 x^2-6x+1& (0,\pm 1), (2,\pm 11) \\
 205 &y^2 =x^6+2 x^4+10 x^3+5x^2-6x+1& (0,\pm 1), (-2,\pm 7)\\
 206 & y^2= x^6 +2 x^4+2 x^3+5 x^2+6x + 1 & (-1, \pm 1), (0, \pm  1), (\frac{1}{2},\pm  \frac{19}{8})\\
 209 & y^2= x^6-4x^5+8x^4-8x^3+8x^2+4x+4& (0,\pm 2),(-\frac{1}{2}, \pm \frac{19}{8})\\
 213 &y^2= x^6+ 2x^4+2 x^3-7x^2+6x-3 &  (1,\pm 1)\\
 215 &y^2=  x^6+4 x^5-12 x^4+20 x^3-20 x^2+12x-4& (1,\pm  1), (2, \pm 10) \\
 221 &y^2= x^6+4x^5+2x^4+6x^3+x^2-2x+1 & (0,\pm  1), (\frac{1}{2}, \pm \frac{9}{8})\\
 230 & y^2=         ( x^3-2 x^2+5x+1) ( x^3+2x^2+x+1)& (0,\pm  1), (1, \pm 5 ), (-2, \pm 5), (3, \pm 35)\\
266 & y^2=
 x^6+4x^5+10x^4+14x^3+17x^2+10x+1& (-1, \pm 1), (0,\pm  1), (-\frac{5}{2},\pm \frac{ 83}{8})\\
285 & y^2=x (x^2+x+4 ) (x^3-x^2-x-3 )& (-1, \pm 4), (0, 0), (3,\pm  24), (-\frac{3}{2}, \pm \frac{57}{8})\\
286 &y^2=(x^3-x^2+3x+1 ) (x^3+x^2-4 ) &(-1,\pm  4), (\frac{5}{2}, \pm\frac{143}{8}) \\
 287 &y^2=  x^6 -4 x^5+2x^4+6x^3-15 x^2+14x-7& (-2,\pm 9) \\
 299 &y^2= x^6-4 x^5+6x^4+6x^3-7x^2-10x-3 & (-\frac{1}{2},\pm \frac{1}{8})\\
 357 & y^2=x^6+8x^4-8x^3+20x^2-12x +12 &(2,\pm 14) \\
 390 &y^2=(x^2-x+1) (x^4+5x^3-8x^2+5x+1 ) &(0, \pm 1), (1,\pm  2)\\ \hline
\end{array}
$$
\begin{centerline}
{Table 1}
\end{centerline}
\end{small}

\subsection{Determination of the rational points}

In order to determine the rational points of the curves $X_0^*(N)$ we will use the so-called elliptic Chabauty method, which uses a
Chabauty procedure on a finite set of unramified $2$-coverings of the curve.

\begin{prop}\label{rational} For the values $N=$ 85, 93, 106, 115, 122, 129, 154, 158, 161, 165, 170, 186, 209, 215, 230, 285, 286,  357
and 390, the set of rational points of $X_0^*(N)(\Q)$ is the set given in Table 1 together with the two points at infinity.
\end{prop}

\noindent {\bf Proof.} The proposition is proved by using some
computations in MAGMA \cite{Magma}. A file with all the computations
can be downloaded from the github account of the third author. We
explain the main ideas in the computation, that can be done for any
hyperelliptic curve $X$ of genus $g$ ($g=2$ in our cases) given by
an equation of the form $y^2=f(x)$, whit $\deg f(x)=2g+2$.

The computation is done in two steps: first one computes the finite
set of twists $C_{\xi}$ of the unramified coverings of the curve $X$
with Galois group $\cong (\Z/2\Z)^{2g}$ which have points locally
for any prime $p$; this is completely analogous to the 2 descent for
elliptic curves, as described in \cite{BrSt09}. Each twists
$C_{\xi}$ is associated to an element $\xi \in (\Q[x]/f(x))^*$
(where the twist corresponding the points at infinity corresponds to
$\xi=1$). If some of the curves does not have rational points
(apparently), one needs to show this by using either a Mordell-Weil
sieve or a higher descent. For our curves this never happens, so we
will not analyze this case further. Our aim now is the determination
of the rational points in $C_{\xi}(\Q)$.

Now, the jacobian of any of this curves has quotients isomorphic to the Weil restriction of elliptic curves $E_{\xi}$ defined over some number fields $K$,
and the rational points in $C_{\xi}(\Q)$ give points in $E_{\xi}(K)$ whose image with respect to a given map
$\varphi_{\xi}:E_{\xi}\to\bP^1$ is in $\bP^1(\Q)$; this is the
necessary data for the elliptic Chabauty function, which computes the set of points in $E_{\xi}(K)$ verifying this condition if
$\rank_{\Z}(E_{\xi}(K))< \deg(K/\Q)$. In practice, the fields $K$ we need are the minimal field of definition of some fixed
factorization $f(x)=g(x)h(x)$ where $g(x)$ has degree $4$.

For example, if $g=2$ and the polynomial $f(x)$ is irreducible, we consider the field $L_0:=\Q[x]/f(x)$. Suppose furthermore that
$f(x)=(x-\alpha)f_1(x)$ in $L[x]$ again with $f_1(x)$ irreducible; then the minimal field of definition of a factorization $f(x)=g(x)h(x)$
where $g(x)$ has degree $4$ is a field of degree $15$ over $\Q$. In this case, the elliptic curves correspond to the jacobians of
the curves $H_{\xi}: y^2=\xi g(x)$, and the map $\varphi_{\xi}: H_{\xi}\to \bP^1$ is given by the
$x$-coordinate. In this case the necessary Chabauty condition is $\rank_{\Z}(E_{\xi}(K))<15$, which is quite likely to be fulfilled;
 but right now the computation of the rank and a finite index subgroup of $E(K)$ for $K/\Q$ of such degree is unfeasible.
 This situation is what happens in all the values of $N$ which we were not able to determine the set of rational points
 (including all prime values of $N$).

The other extreme case is when there is a factorization $f(x)=g(x)h(x)$ already defined over $\Q$; in this case we need to compute the rational points
of the curves $y^2=d_{\xi} g(x)$ for some values $d_\xi\in\Q^*$, which are only finite if the rank is zero (which is very unlikely to happen for all
the necessary twists $d_\xi$).

The best cases are when one can find such a factorization over a
field of degree $\le 4$ for any twists verifying the corresponding
Chabauty condition. In some cases we used distinct fields for
different twists, as we explain in the following example.

\begin{exemp} We explain in detail the case $X_0^*(85)$, where $f(x)=(x^2 - 2x + 5)(x^4 - 2x^3 + 3x^2 - 6x + 5)$. We have
$$X_0^*(85)(\Q)=\{(0,\pm 5),(1,\pm 2), (2,\pm 5), (\frac{3}{2},\pm \frac{17}{8}), (-\frac{4}{3},\pm \frac{425}{27}),\pm \infty \}.$$
We have 5 twists, corresponding to the $x$-coordinates of the
rational points, except that the points with $x=1$ already appear in
the trivial twists corresponding to the points at infinity.

If we consider the given factorization over $\Q$, the corresponding
elliptic curves for all the twists except the trivial one have rank
1. On the other hand the curve $H$ given by the equation $y^2=x^4 -
2x^3 + 3x^2 - 6x + 5$ has rank 0 and 4 points, corresponding to the
points at infinity and the points with $x=1$.

If we adjoint a root of $x^2 - 2x + 5$ we get a quadratic extension $K_2/\Q$. Over this extension we get also a factorization
of  $x^4 - 2x^3 + 3x^2 - 6x + 5=h_1(x)h_2(x)$ with $\deg(h_i(x))=2$ for $i=1,2$. So we can take $f(x)=g(x)h_1(x)$ for some $g(x)$ of degree 4.
The twists corresponding to the points with $x$-coordinate $\frac{3}{2}$ and $-\frac{4}{3}$ have rank 1, and Chabauty method succeeds.
But the ones corresponding to the points with $x$-coordinate $0$ and $2$ have rank 2.

If we adjoint a root of $x^4 - 2x^3 + 3x^2 - 6x + 5$ we get a field
$K_4$ where $f(x)$ has $4$ roots and a degree $2$ factor. By
considering the corresponding degree 4 polynomial as a product of
two (adequate) degree one factors and the degree two factor we
finally get that the remaining twist have jacobian of rank 1 and we
find a non-torsion point in each case, and Chabauty computations
succeeds.
\end{exemp}

There is one case were the approach described above did not succeed.

\begin{exemp} In the case $X_0^*(390)$ we had to do a slightly modified approach; in fact, we tried fields of
degree 1 and 2 and the Chabauty condition was not fulfilled, and we had to go to a degree 8 extension, where the rank computation did not succeed.

Instead, we showed that all the rational points in $X_0^*(390): y^2=(x^2-x+1)(x^4+5x^3-8x^2+5x+1)$ came from the trivial twists over $\Q$. This means
that their $x$-coordinates verify that $$  y_1^2=x^2-x+1 \text{ and } y_2^2=x^4+5x^3-8x^2+5x+1$$ for some $y_1,y_2\in \Q$. The curve $X'$ determined by
these equations is an hyperelliptic curve of genus $3$, whose hyperelliptic equation can be computed by parametrizing the first equation. We get
the equation $$X': z^2=t^8 + 10t^7 - 41t^6 + 42t^5 + 33t^4 - 76t^3 + 44t^2 - 16t + 4.$$

Now we apply the above method to this new curve $X'$. The curve has
(apparently) $12$ rational points (two above each rational point of
$X_0^*(390)$, as it is an unramified 2-covering), with $6$ possible
values for the $x$-coordinates. We computed there are exactly 6
possible twists, one for each $x$-coordinate. Over $K=\Q(\sqrt{5})$
the defining hyperelliptic polynomial factors as a product of two
degree 4 polynomials. For every twists one of the two quotient
elliptic curves of the corresponding covering has rank one, and the
Chabauty computation succeeds. \hfill $\Box$
\end{exemp}

\section{An example and results}
First, we show, for the case $N=67$, the procedure used.
The equation obtained in Table 1  is $y^2=x^6 - 4 x^5 + 6 x^4 - 6 x^3 + 9 x^2 - 14 x + 9 $. Let $j_{67}(q):=j(q^{67})=1/q^{67}+ 744+\cdots$ and put $J_1=j+j_{67}$ and $J_2=j\cdot j_{67}$. With the notation of before section, we have
$$
f_3=\displaystyle{\frac{1}{2} \left (-1 + x - 2 x^2 + x^3 + y\right )}\,,\quad f_4=x\,f_3+1\,,\quad f_5=x\,f_4-1\,.
$$
After computing, we obtain
\begin{scriptsize}
$$\begin{array}{cr}
J_1=&-23 f_3^{22} + f_3^{21} f_4 - 1279 f_3^{21} - 781 f_3^{20} f_4 + 186 f_3^{20} f_5 - 99914 f_3^{20} -
 39399 f_3^{19} f_4 + 14954 f_3^{19} f_5 - 2698696 f_3^{19} - 265633 f_3^{18} f_4 +  \\ [3 pt]
& 380472 f_3^{18} f_5 - 20514523 f_3^{18} + 6929641 f_3^{17} f_4 + 1576893 f_3^{17} f_5 - 49240824 f_3^{17} +
  67627402 f_3^{16 }f_4 - 16450546 f_3^{16} f_5 - \\ [3 pt]
&61401116 f_3^{16} +
 190686364 f_3^{15} f_4 - 81315034 f_3^{15} f_5 - 56264079 f_3^{15} + 259977664 f_3^{14} f_4 -  148558638 f_3^{14 }f_5 - 88533538 f_3^{14 }+\\[ 3 pt] &
 95806265 f_3^{13 }f_4 - 69608123 f_3^{13} f_5 -  162557463 f_3^{13} -
 295479289 f_3^{12} f_4 + 158123161 f_3^{12} f_5 - 27169544 f_3^{12} -
 558873206 f_3^{11} f_4 +\\ [ 3 pt] & 260674425 f_3^{11} f_5 + 456803156 f_3^{11} -
 423722114 f_3^{10} f_4 + 202709065 f_3^{10} f_5 + 731171796 f_3^{10} -
 51627779 f_3^9 f_4 + 133780373 f_3^9 f_5 +\\[ 3 pt] & 234273070 f_3^9 +
 264268555 f_3^8 f_4 - 64460559 f_3^8 f_5 - 502745764 f_3^8 +
 337312727 f_3^7 f_4 - 318069668 f_3^7 f_5 - 623447279 f_3^7 +\\[ 3 pt] &
 158991342 f_3^6 f_4 - 299948399 f_3^6 f_5 - 229889114 f_3^6 -
 28504966 f_3^5 f_4 - 91453878 f_3^5 f_5 + 60433254 f_3^5 -
 60041832 f_3^4 f_4 + \\[3 pt] &24362830 f_3^4 f_5 + 79628320 f_3^4 -
 20570848 f_3^3 f_4 + 26593344 f_3^3 f_5 + 23436576 f_3^3 -
 941600 f_3^2 f_4 + 8600928 f_3^2 f_5 + 1047456 f_3^2 + \\ [ 3 pt] &647200 f_3 f_4 +
 1386464 f_3 f_5 - 571936 f_3 + 81536 f_4 + 92000 f_5 - 65536\,,
 \end{array}
$$
\end{scriptsize}
and
\begin{scriptsize}
$$\begin{array}{cr}
J_2=&f_5f_3^{21} + 720 f_4f_3^{21} + 179980f_3^{22} + 17300122 f_5f_3^{20} +
 410510311 f_4f_3^{20} + 5149868567f_3^{21} + 42380978353 f_5f_3^{19} +\\ [ 3 pt] &
 176848900626 f_4f_3^{19} + 839384847849f_3^{20 }+
 3347632163474 f_5f_3^{18} + 5490095012794 f_4f_3^{18} +
 20232433296285f_3^{19} + \\ [ 3 pt] & 58903607428273 f_5f_3^{17} +
 52642797600751 f_4f_3^{17} + 194740838424278f_3^{18} +
 465038243745693 f_5f_3^{16} + \\ [ 3 pt] &252244831282013 f_4f_3^{16} +
 1049935814900775f_3^{17} + 2137932776610224 f_5f_3^{15} +
 718160756825707 f_4f_3^{15} + \\ [ 3 pt] &3682889678423580f_3^{16} +
 6483265127099687 f_5f_3^{14 }+ 1296256091753552 f_4f_3^{14} +
 9112585010431660f_3^{15} + \\ [ 3 pt] &13905322982585428 f_5f_3^{13} +
 1455284181206971 f_4f_3^{13} + 16698439687036025f_3^{14} +
 22013212259613947 f_5f_3^{12} +\\ [ 3 pt]  & 785970239416766 f_4f_3^{12} +
 23380324317471236f_3^{13 } + 26434965501463926 f_5f_3^{11 }-
 389286347179143 f_4f_3^{11} + \\ [ 3 pt] &25533239990035759f_3^{12} +
 24511794195150313 f_5f_3^{10} - 1257772458271356 f_4f_3^{10 }+
 22044986048091002f_3^{11 }+ \\ [ 3 pt] &17745745315071401 f_5f_3^9
 1369085873848977 f_4f_3^9 + 15174609069161127f_3^{10} +
 10090414014987393 f_5f_3^8 -\\ [ 3 pt] & 954179629348608 f_4f_3^8 +
 8365619455927489f_3^9 + 4511870660986624 f_5f_3^7 -
 476646738132432 f_4f_3^7 + 3699701950122944f_3^8 +\\ [ 3 pt] &
 1580542220676176 f_5f_3^6 - 178202666014272 f_4f_3^6 +
 1312789289946672f_3^7 + 429477284579616 f_5f_3^5 -
 51157416030048 f_4f_3^5 + \\ [ 3 pt] &374258702810464f_3^6 +
 88841374862944 f_5f_3^4 - 11506808706816 f_4f_3^4 +
 86310527850080f_3^5 + 13536732334080 f_5f_3^3 -\\ [ 3 pt] &
 2038076477952 f_4f_3^3 + 16327740526336f_3^4 +
 1432454117376 f_5f_3^2 - 272275906560 f_4f_3^2 + 2547759661312f_3^3 +\\ [ 3 pt] &
 93842541824 f_5f_3 - 24005742848 f_4f_3 + 313280547584f_3^2 +
 2852000000 f_5 - 1009741824 f_4 + 26269884672f_3 + 1073741824\,.
 \end{array}
$$
\end{scriptsize}
Hence, the $j$-invariants of the $\Q$-curves attached to $\infty'$ are the solutions of the equation
$$
z^2-J_1(\infty') \,z+J_2(\infty')=0\,.
$$
Since $(J_1(\infty'),J_2(\infty'))=(- 65536,1073741824)$, we get $j=-32^3$, which corresponds to an elliptic curve with CM  by the quadratic order of discriminant $-11$.
The remaining rational non-cuspidal points provide $\Q$-curves with CM:
$$
\begin{array}{c|r||c|r||c|r||c|r||}
\text{point}& j &\text{point}& j&\text{point}& j &\text{point}& j \\\hline\hline
(-1,\phantom{-}7) & 255^3 &(0,\phantom{-}3)& -3\cdot 160^3 & (1,\phantom{-}1)& 20^3 & (2,\phantom{-}1) &-960^3\\
(-1,-7)&-5280^3 &(0,-3) &0 & (1,-1)&-15^3 &(2,-1) &2\cdot 30^3\\ \hline\hline
\end{array}
$$
Next, we show the results obtained.
\subsection{$N$ is a  prime: $\Gal (\Q( j)/\Q)\hookrightarrow\Z/2\Z$.}
$$
\begin{array}{c|c|c|r|r|}
N&\text{ point} & \text{CM}&D&j\text{-invariant}\phantom{ccccc}\\ \hline\hline
67 & \infty' &\text{yes}& -11 &-32^3\\
&(-1,7)&\text{yes}& -7&255^3 \\
&(-1,-7)&\text{yes} &-67&-5280^3\\
&(0,3)&\text{yes} &-27 &-3\cdot 160^3\\
&(0,-3)&\text{yes} &-3 &0 \\
&(1,1)&\text{yes} & -8 &20^3 \\
&(1,-1)&\text{yes} &-7& -15^3 \\
&(2,1)&\text{yes} &-43& -960^3 \\
&(2,-1) &\text{yes}&-12 &2\cdot 30^3 \\ \hline\hline
73 & \infty' &\text{yes}&-12 &2\cdot 30^3 \\
&(0,1)&\text{yes}&-27 &-3 \cdot 160^3 \\
&(0,-1&\text{yes})&-4 &12^3 \\
&(1,1)&\text{yes} &-19&- 96^3\\
&(1,-1)&\text{yes} &-8&20^3 \\
&(2,3) &\text{yes}&-67& -5280^3 \\
&(2,-3)&\text{yes} &-16& 66^3 \\
&(3/2,5/8)&\text{yes} & -3&0 \\
&(3/2,-5/8)&\text{non}&- &\displaystyle{20 \left(\frac{ 3 (-26670989 \pm  15471309 \sqrt{-127})}{2^{26}}\right)^3}\\ \hline\hline
\end{array}
$$

$$
\begin{array}{c|c|c|r|r|}
N&\text{ point} & \text{CM}&D&j\text{-invariant}\phantom{ccccc}\\ \hline\hline
103 & \infty'&\text{yes} &-67&-5280^3 \\
&(0,1)&\text{yes}& -43&-960^3 \\
&(0,-1)&\text{yes}&-27&-3\cdot 160^3 \\
&(1,1) &\text{yes}&-19&-96^3 ]\\
&(1,-1)&\text{yes} &-12&2\cdot 30^3  \\
&(3,19)&\text{non} &-&19 (48 (1623826405 \pm 30228849 \sqrt{2885}))^3  \\
&(3,-19)&\text{yes} &-3& 0 \\ \hline\hline
107 & \infty' &\text{yes} & -8& 20^3 \\
  & (0,1) &\text{yes} & -7& -15^3\\
 & (0,-1) &\text{yes} & -43&-960^3 \\
  & (2,1) &\text{yes} &-67& -5280^3 \\
 & (2,-1) &\text{yes} &-28& 255^3 \\  \hline\hline
167 & \infty'&\text{yes}  &-43&  -960^3  \\
 & (-1,1)&\text{yes} & -67&-5280^3  \\
 &(-1,-1) &\text{yes}&-163&-640320^3 \\ \hline\hline
191 & \infty' &\text{yes} &-43& -960^3 \\
 &(0,1) &\text{yes}&-11&-32^3 \\
 &(0,-1)&\text{yes} &-7&-15^3\\
  &(2,11) &\text{non}& -&j_0 \\
 &(2,-11)&\text{yes} &-28&255^3\\
\hline\hline
\end{array}
$$
where
\begin{scriptsize}
$$j_0=( 724537954586714121 \pm16056976492100\sqrt{2036079533})\left(\frac{480(7725788647437 \pm 95942438\sqrt{2036079533})}{191^2} \right)^3
$$
\end{scriptsize}
\subsection{ $N$ is a product of two primes: $\Gal (\Q( j)/\Q)\hookrightarrow (\Z/2\Z)^2$.}
$$
\begin{array}{c|c|r|r|r|}
N&\text{ point} & \text{CM }&D& $j$ \text{ or }\Q(j)\\ \hline\hline
85 & \infty' &\text{yes} &-19 &-96^3\\
&(0,5)  &\text{yes}&-35 & -(16(15 \pm  7 \sqrt 5))^3 \\
&(0,-5)  &\text{yes}&-60 & (3(470 \pm  213 \sqrt 5)^3 (1 \pm  \sqrt 5)/2)\\
&(1,2) &\text{yes}&-16 & 66^3 \\
&(1,-2)  &\text{yes}& -4 &12^3 \\
&(2,5)  &\text{yes}&-115 & -(48 (785 \pm  351 \sqrt 5))^3 \\
&(2,-5) &\text{yes} & -15&-(3 (25 \pm 9 \sqrt 5)/2)^3 (-1 \pm \sqrt 5)/2\\
&(3/2,17/8)  &\text{yes}&-51& -(48 (37 \pm 9 \sqrt{17}) )^3 (-4 \pm \sqrt {17}) \\
&(3/2,-17/8) &\text{non} & - &  \Q(\sqrt{17},\sqrt{-95})\\
&(-4/3,425/27)  &\text{non}& - &  \Q(\sqrt{85},\sqrt{-4295})\\
&(-4/3,-425/27) &\text{yes} &-595 & \Q(\sqrt5,\sqrt {17})\\
\hline\hline
\end{array}
$$

$$
\begin{array}{c|c|r|r|r|}
N&\text{ point} & \text{CM }&D& $j$ \text{ or }\Q(j)\\ \hline\hline
93 & \infty'  &\text{yes}& -12 &2\cdot 30^3 \\
 & (0,3) &\text{yes}&-60 &(3(470 \pm  213 \sqrt 5)^3 (1 \pm  \sqrt 5)/2)\\
&(0,-3)  &\text{yes}&-24 &(12 (5 \pm  2 \sqrt 2))^3 (3 \pm  2 \sqrt 2)\\
& (-1,3)  &\text{yes}&-123 &-(480 (461 \pm  72 \sqrt{41}))^3 (-32 \pm  5  \sqrt{41})
\\
&(-1,-3) &\text{yes} &  -75 &-(48 (-69 \pm  31 \sqrt 5))^3 (\pm \sqrt 5)\\
 &(1,1)  &\text{yes}&-11 &-32^3 \\
  & (1,-1) &\text{yes} &-3,-12& 0, -3\cdot 160^3  \\
   &(2,3) &\text{yes} &-147 &-3 (480 (142 \pm  31 \sqrt{21}))^3(\pm \sqrt{21})
 \\
   &(2,-3) &\text{yes} &-15 &-(3 (-5\pm 4 \sqrt 5 ))^3 (-3\pm\sqrt 5 )/2\\
   & (3/2,9/8)  &\text{yes}&-48&4 (15 (30 \pm 17 \sqrt 3))^3 \\
    &(3/2,-9/8)  &\text{non}&-  & \Q(\sqrt{-15},\sqrt{-109}) \\
     &(1/4,143/64) &\text{yes}&-3 & 0 \\
     &(1/4,-143/64) &\text{non}&- & \Q(\sqrt{-23},\sqrt{-143})  \\ \hline\hline
     106 & \infty'&\text{yes} &-7 & -15^3  \\
 &(-1,4) &\text{yes}&-36 &4 (21\pm 20 \sqrt3 )^3 (7 \pm 4 \sqrt 3)\\
  &(-1,-4)&\text{yes}&-148 &(60 (2837 \pm  468 \sqrt{37}))^3\\
   & (0,1)&\text{yes}&-40 &(6(65 \pm  27 \sqrt 5) )^3   \\
    &(0,-1) &\text{yes}&-4,-16 &12^3, 66^3  \\
     & (1,2)&\text{yes} &-24 & (12 (9 \pm 7 \sqrt 2))^3 (-1 \pm \sqrt 2)  \\
     &(1,-2)&\text{yes} &-52 &(30 (31 \pm  9 \sqrt {13}))^3 \\
      &(2,5)&\text{yes}&-100 &(6 (2927 \pm  1323 \sqrt 5))^3 \\
      & (2,-5)&\text{yes} &-4 &12^3  \\
       &(1/2,5/8)&\text{non} & - & \Q(\sqrt{33},\sqrt{-591})\\
        &(1/2,-5/8)&\text{yes}&-7, -28 &-15^3,255^3\\
  \hline\hline
  115 & \infty' &\text{yes}&-115 & (48 (-785 \pm  351 \sqrt 5))^3\\
 &(1,1)&\text{yes}&-19& -96^3\\
  &(1,-1)&\text{yes}&-11 & -32^3\\
   & (2,5)&\text{yes}& -235 &(528 (-8875 \pm  3969 \sqrt 5))^3  \\
    &(2,-5)&\text{yes} &-15 &-( 3/2 (25 \pm 9 \sqrt 5))^3 (-1 \pm \sqrt 5)/2\\
     & (1/2,5/8) &\text{non}&- &\Q(\sqrt{65},\sqrt{-3495})\\
     &(1/2,-5/8)&\text{yes} &-40 &(6 (65 \pm  27 \sqrt 5))^3 \\
     & (4/3,35/27)&\text{yes} &-60 &(3(470 \pm  213 \sqrt 5)^3 (1 \pm  \sqrt 5)/2)
      \\
     &(4/3,-35/27) &\text{non}&-&  \Q (\sqrt{10},\sqrt{-9278})\\
     \hline\hline
      122 & \infty' &\text{yes}&-36 &- (4 (102 \pm 61 \sqrt 3))^3 (-2 \pm \sqrt 3)  \\
     &(-1,4)&\text{yes}&-52 &(30 (31 \pm  9 \sqrt{13}))^3\\
 &(-1,-4)&\text{yes}&-100 &(6 (2927 \pm  1323 \sqrt 5))^3\\
   & (0,1)&\text{yes}&-3,-12 & 0, 2\cdot 30^3  \\
    &(0,-1)&\text{yes} &-4,-16 &12^3, 66^3  \\
     & (1,2)&\text{yes} &-88&(60 (155 \pm  108 \sqrt 2))^3 \\
     &(1,-2)&\text{yes} &-20 &(2 (25 \pm  13 \sqrt 5) )^3 \\
     & (3/2,37/8)&\text{yes} & -232 &(30 (140989 \pm  26163 \sqrt {29}))^3 \\
     &(3/2,-37/8)&\text{non} &-&\Q(\sqrt{-15},\sqrt{ 1585}) \\
       & (2/3,37/27)&\text{yes} & -4 &12^3   \\
     &(2/3,-37/27)&\text{non} &-&\Q(\sqrt{1258}, \sqrt{-1598})   \\
     \hline\hline
\end{array}
$$

$$
 \begin{array}{c|c|r|r|r|}
N&\text{ point} & \text{CM}&D & j \text{ or }\Q(j)\\ \hline\hline
      129 & \infty'&\text{yes} &-75 &-(48 (69 \pm 31 \sqrt 5))^3(\pm  \sqrt 5 )   \\
 &(-1,3)&\text{yes}&-123&-(480 (-461 \pm 72 \sqrt{41}))^3 (32 \pm  5 \sqrt{41})\\
 &(-1,-3)&\text{yes}&-48 &4 (15 (30 \pm 17 \sqrt 3 ))^3    \\
   & (0,2)&\text{yes}&-147 &-3 (480 (362 \pm 79 \sqrt{21}))^3(\pm  \sqrt{21})
   \\
    &(0,-2) &\text{yes}&-8 &20^3  \\
     & (1,1) &\text{yes}&-3,-27 &0, -3\cdot 160^3 \\
          &(1,-1)&\text{yes} &-12 &2\cdot 30^3 \\
     & (1/2,3/8)&\text{yes} &-51&-(48 (37 \pm  9 \sqrt{17}))^3 (-4 \pm \sqrt{17})   \\
     &(1/2,-3/8) &\text{non}&-& \Q(\sqrt{57}, \sqrt{-687})  \\
       &  (-7/5, 383/125)&\text{yes}&-3 &0  \\
     &(-7/5, -383/125) &\text{non}&-& \Q(\sqrt{1149},\sqrt{-1059})   \\
       & (7/12, 383/1728)&\text{non} & - & \Q(\sqrt{-7},\sqrt{-444783}) \\
     & (7/12, -383/1728)&\text{non}&-&\Q(\sqrt{85},\sqrt{-347 })  \\
     \hline\hline
     133 & \infty' &\text{non}& -& \Q (\sqrt 2,\sqrt{69}) \\
 &(0,1)&\text{yes}&-27 & -3\cdot 160^3 \\
  &(0,-1)&\text{yes}&-19& -96^3\\
   & (1,1))&\text{yes}&-91 &(48 (-227 \pm  63 \sqrt{13}))^3   \\
    &(1,-1)&\text{yes} &-12 &2\cdot 30^3 \\
     & (3/5, 83/125)&\text{yes} &-3 & 0 \\
     &(3/5, -83/125) &\text{non}&-&   \Q(\sqrt{-31},\sqrt{-3651})\\
          \hline\hline
134 & \infty'&\text{yes} &-52 &(30 (31 \pm  9 \sqrt{13}))^3  \\
 &(-1,3)&\text{yes}& -7&-15^3 \\
  &(-1,-3)&\text{yes}& -232 &(30 (140989 \pm  26163 \sqrt{29}))^3\\
   & (0,1))&\text{yes}&-20&(2 (25 \pm  13 \sqrt 5))^3   \\
    &(0,-1) &\text{yes}&-3, -12 &0, 2\cdot 30^3 \\
     & (1, 1)&\text{yes} &  -8&20^3\\
     &(1,-1) &\text{yes}&-7,-28 &-15^2, 255^3   \\
     & (-1/2, 7/8)&\text{non} &- & \Q(\sqrt{113},\sqrt{-1271}) \\
     &(-1/2, -7/8)&\text{yes} &-72&(20 (389 \pm  158 \sqrt 6))^3 (-5 \pm 2 \sqrt 6)   \\
          \hline\hline
  146 & \infty' &\text{yes}&-3,-12&0, 2\cdot 30^3  \\
 &(-1,1)&\text{yes}&-36& -(4 (102 \pm  61 \sqrt 3))^3 (-2 \pm  \sqrt 3) \\
  &(-1,-1)&\text{yes}& -148&(60 (2837 \pm 468 \sqrt {37}))^3\\
   & (0,1))&\text{yes}&-4, -16 &12^3,66^3  \\
    &(0,-1)&\text{yes} &-24 &(12 (9 \pm  7 \sqrt 2))^3 (-1 \pm  \sqrt 2) \\
     & (1, 3) &\text{yes}& -8&20^3 \\
     &(1,-3)&\text{yes} &-72&-(20 (389 \pm  158 \sqrt 6))^3 (-5 \pm 2 \sqrt 6)    \\
     & (2, 5)&\text{yes} &-100&(6 (2927 \pm  1323 \sqrt 5))^3  \\
     &(2, -5)&\text{yes} &-4 &12^3   \\
          \hline
 \hline
 158 & \infty' &\text{yes}&-7 & -15^3  \\
& (0, 1)&\text{yes} &-3,-12& 0, 2\cdot 30^3   \\
 &(0,-1)&\text{yes}& -24 &(12 (9 \pm  7 \sqrt 2))^3 (-1 \pm  \sqrt 2) \\
  &    (2,1))&\text{yes}& -232 &(30 (140989 \pm  26163 \sqrt{29}))^3   \\
    &(2,-1)&\text{yes} &-148 &(60 (2837 \pm  468 \sqrt{37}))^3 \\
     & (1/2, 1/8)&\text{non} & - &\Q(\sqrt{1169},\sqrt{-1247})  \\
     &(1/2,-1/8)&\text{yes} &-7,-28 &-15^3, 255^3   \\
               \hline      \hline
     \end{array}
     $$

$$
 \begin{array}{c|c|r|r|r|}
N&\text{ point} & \text{CM}&D & j \text{ or }\Q(j)\\ \hline\hline
           161 & \infty' &\text{yes}&-7&-15^3   \\
 &(-1,7)&\text{yes}& -91&-(48 (227 \pm  63 \sqrt{13}))^3  \\
  &(-1,-7)&\text{yes}& -483& \Q(\sqrt{21},\sqrt{69})\\
   & (1,1))&\text{yes}&-115& -(48 (785 \pm 351 \sqrt 5))^3  \\
    &(1,-1)&\text{yes} &-19 &-96^3\\
     & (-1/2, 35/8) &\text{non}&&\Q(\sqrt{-7},\sqrt{32009}) \\
     &( -1/2, -35/8)&\text{yes}&-112 &(15 (2168 \pm 819 \sqrt 7))^3   \\
     & (-1/4, 209/64)&\text{yes} & -8&255^3 \\
     & (-1/4, -209/64)&\text{non}&&\Q)(\sqrt{209},\sqrt{  -1140391}) \\
          \hline\hline
          177 & \infty' &\text{yes}&--11 &32^3  \\
 &(0,1)&\text{yes}&-24&(12 (9 \pm  7 \sqrt 2))^3 (-1 \pm  \sqrt 2)  \\
  &(0,-1)&\text{yes}&-8 & 20^3 \\
   & (3/2, 17/8)&\text{yes}&-267&-(240 (562501 \pm  59625 \sqrt{89}))^3 (-500 \pm 53 \sqrt{89})   \\
    &(3/2, -17/8)&\text{non} &-&\Q(\sqrt{-23},\sqrt{2881})  \\
               \hline\hline
               205 & \infty' &\text{yes}&-115 &-(48 (785 \pm  351 \sqrt 5))^3  \\
 &(0,1)&\text{yes}&-16 &66^3  \\
  &(0,-1)&\text{yes}&-40 & (6 (65 \pm  27 \sqrt 5))^3 \\
   & (-2,7)&\text{yes}& -4&12^3   \\
    &(-2,-7) &\text{yes}&-1435& \Q(\sqrt 5 ,\sqrt{21})\\
               \hline\hline
    206 & \infty' &\text{yes}&-24&(12 (9 \pm  7 \sqrt 2))^3 (-1 \pm  \sqrt 2)  \\
 &(-1,1)&\text{yes}& -3,-12&0,2\cdot 30^3 \\
  &(-1,-1)&\text{yes}&-88&(60 (155 \pm  108 \sqrt 2))^3  \\
   & (0,1)&\text{yes}&-40& 6 (65 \pm  27 \sqrt 5 )    \\
    &(0,-1)&\text{yes} &-20&(2 (25 \pm  13 \sqrt 5))^3  \\
      & (1/2,19/8)&\text{yes}&-148& (60 (2837 \pm  468 \sqrt {37}))^3   \\
    &(1/2,-19/8)&\text{non} &-&\Q (\sqrt{193},\sqrt{-27119})\\
               \hline   \hline
      209 & \infty'&\text{yes} &-8&20^3   \\
 &(0,2)&\text{yes}& -19 &-96^3  \\
  &(0,-2)&\text{yes}&-88 &(60 (155 \pm 108 \sqrt 2))^3 \\
      &(-1/2,19/8)&\text{non}    &- & \Q(\sqrt{-1007},\sqrt{902537}) \\
    &(-1/2,-19/8)&\text{yes} &-627&\Q( \sqrt{33},\sqrt{57})\\
               \hline           \hline
               213 & \infty' &-51&\text{yes}&-(48 (37 \pm  9 \sqrt{17}))^3 (-4 \pm  \sqrt{17}) \\
 &(1,1)&\text{yes}&-123&(480 (461 \pm  72 \sqrt{41}))^3 (-32 \pm  5\sqrt{41})  \\
  &(1,-1)&\text{yes}&-11& -32^3  \\
\hline
\hline
215 & \infty'&\text{non} &- &\Q (\sqrt 2,\sqrt{47645})      \\
 &(1,1)&\text{yes}&-235 &-(528 (8875 \pm  3969 \sqrt 5))^3 \\
  &(1,-1)&\text{yes}& -19 &-96^3\\
      &(2,10)&\text{non}    &- & \Q(\sqrt{85},\sqrt{3418805})  \\
    &(2,-10)&\text{yes} &-115 &-(48 (785 \pm  351 \sqrt 5))^3 \\
               \hline\hline
             221 & \infty'&\text{yes}&-16&66^3  \\
 &(0,1)&\text{yes}&-43&-960^3   \\
 &(0,-1)&\text{yes}&-51&-(48 (37 \pm 9 \sqrt{17}))^3 (-4 \pm  \sqrt{17})  \\
       &(1/2,9/8)&\text{non}    &- & \Q(\sqrt{1081},\sqrt{ -779263} )  \\
    &(1/2,-9/8)&\text{yes} &-4&12^3  \\
               \hline  \hline
 \end{array}
$$

  $$
 \begin{array}{c|c|r|r|r|}
N&\text{ point} & \text{CM}&D & j \text{ or }\Q(j)\\ \hline\hline
   287 & \infty'&\text{yes}&-91 &-(48 (227 \pm 63 \sqrt{13}))^3   \\
 &(-2,9)&\text{yes}& -1435& \Q(\sqrt 5,\sqrt{41}  )  \\
 &(-2,-9)&\text{non}&-&\Q(\sqrt{8321},\sqrt{2904137173}) \\
 \hline\hline
 299 & \infty'&\text{yes}&-91 &-(48 (227 \pm 63 \sqrt{13}))^3    \\
 &(-1/2,1/8)&\text{yes}& -43&-960^3    \\
 &(-1/2,-1/8)&\text{non}&&\Q( \sqrt{1513},\sqrt{-3325543}) \\ \hline\hline
     \end{array}
     $$

\subsection{$N$ is a product of three primes: $\Gal (\Q( j)/\Q)\hookrightarrow(\Z/2\Z)^3$.}
\begin{small}
$$\begin{array}{c|c|c|c|c|}
N&  \text{point} & \text{CM}  & D &  j\text{ or }\Q(j)\\ \hline\hline
154 &\infty' &  \text{yes}  & -40&(6 (65 \pm 27 \sqrt 5))^3\\
 &(0,2)&   \text{yes}       & -24& (12 (9 \pm  7 \sqrt 2))^3 (-1 \pm  \sqrt 2)\\
&(0,-2)&    \text{yes}      & -52&(30 (31 \pm  9 \sqrt{13}))^3\\
&(1,4)     & \text{yes}   & -7&-15^3\\
&(1,-4)     &    \text{yes}&-7,-28 & -15^3,255^3\\
&(2,0) &  \text{yes}&-84 & \Q(\sqrt 3,\sqrt 7)\\
&(-3/2,77/8) &\text{non}&--&\Q(\sqrt{-143},\sqrt{-185},\sqrt{-455})\\
&(-3/2,-77/8)&\text{yes}&-1540 &\Q(\sqrt5,\sqrt 7,\sqrt{11}) \\
&(-1/3,56/27) &\text{non}&-- &\Q(\sqrt{7},\sqrt{5\cdot 11},\sqrt{-479})\\
 &(-1/3,-56/27)&\text{yes}&-28,-112  &255^3,(15 (2168 \pm 819 \sqrt 7))^3\\
& (4,22) &\text{yes}&-1848 &\Q(\sqrt 2,\sqrt{21},\sqrt{33})\\
 &(4,-22) &\text{yes}& -132&\Q(\sqrt 3,\sqrt{11})\\ \hline\hline
 165 &\infty'& \text{yes}   &  -11   &-32^3\\
    &(0,  3) &\text{yes}&-195&\Q(\sqrt 5,\sqrt{13})\\
  &(0,  -3) &\text{yes}&-51&\-(48 (37 \pm 9 \sqrt{17}) )^3 (-4 \pm \sqrt {17})\\
  &(1,0) &\text{yes}&-24&(12 (9 \pm  7 \sqrt 2))^3 (-1 \pm \sqrt 2)\\
  &(2,  5) &\text{yes}&-435 &\Q(\sqrt 5,\sqrt{29})\\
  &(2,  -5) &\text{yes}&-35&-(16 (15 \pm  7 \sqrt 5))^3\\
  &(-1/2,  15/8) &\text{non}&--&\Q(\sqrt{-15},\sqrt{265},\sqrt{1745})\\
  &(-1/2,  -15/8) &\text{yes}&-120 & \Q(\sqrt 2,\sqrt 5)\\
     &(-3,0) &\text{yes}&-1155&\Q (\sqrt 5,\sqrt{21},\sqrt{33})\\
      &(2/3,55/27) &\text{yes}&-11,-99&-32^3, (16 (3751 \pm 653 \sqrt{33}))^3 (-23 \pm 4 \sqrt{33})\\
      &(2/3,-{55}/{27}) &\text{non}&--&\Q(\sqrt{-11},\sqrt{47},\sqrt{-661})\\
      &(5/2,99/8 &\text{yes}&-1320&\Q(\sqrt5,\sqrt 6,\sqrt{22})\\
      &(5/2,-99/8) &\text{non}&--& \Q(\sqrt{-7},\sqrt{33},\sqrt{393})\\ \hline \hline
170 &\infty'&\text{yes}&-36& -(4 (102 \pm  61 \sqrt 3))^3 (-2 \pm  \sqrt 3)\\
&(-1, 2) &\text{yes}&-4,-16& 12^3,66^3\\
&(-1, -2) &\text{yes}&-340&\Q(\sqrt 5,\sqrt{17})\\
&(0, 1) &\text{yes}&-4,-100&12^3, (6 (2927 \pm 1323 \sqrt 5))^3\\
&(0, -1) &\text{yes}&-15&-(3 (25 \pm 9 \sqrt 5)/2 )^3 (-1 \pm \sqrt 5)/2\\
 &(2, 5) &\text{yes}&-280&\Q(\sqrt 2,\sqrt 5)\\
 & (2, -5) &\text{yes}&-15,-60&-(3 (25 \pm 9 \sqrt 5)/2 )^3 (-1 \pm \sqrt 5)/2,-(3 (470 \pm 213 \sqrt 5))^3(1 \pm \sqrt 5)/2 \\
 &(-1/2, 5/8)&\text{yes}&-120&\Q(\sqrt2,\sqrt 5)\\
 & (-1/2, -5/8) &\text{non}&--&\Q(\sqrt{17},\sqrt{-95},\sqrt{65})\\
 &(5/3, 38/27) &\text{non}&--&\Q(\sqrt{73},\sqrt{19},\sqrt{-5})\\
 & (5/3, -38/27) &\text{yes}&-4&12^3\\ \hline\hline
\end{array}
$$
\newpage
$$\begin{array}{c|c|c|c|c|}
N&  \text{point} & \text{CM}  & D & j \text{ or }\Q(j)\\ \hline\hline
186 &\infty' &  \text{yes}  & -3,-12 &0, 2\cdot 30^3\\
&(-1,-3)&    \text{yes}      &-228&\Q(\sqrt 3,\sqrt {19})\\
&(0,1)     & \text{yes}        &-15&-(3 (25 \pm 9 \sqrt 5)/2 )^3 (-1 \pm \sqrt 5)/2\\
&(0,-1)     &    \text{yes}      &-24&(12 (9 \pm  7 \sqrt 2))^3 (-1 \pm  \sqrt 2)\\
&(1,3) &  \text{yes}      &-168&\Q(\sqrt 6,\sqrt{14})\\
&(1,-3) &\text{yes}&-120&\Q(\sqrt 2 ,\sqrt 5)\\
&(2,9) &\text{yes}& -708&\Q(\sqrt 3,\sqrt {59})\\
&(2,-9)&\text{yes}       &-15,-60&-(3 (25 \pm 9 \sqrt 5)/2 )^3 (-1 \pm \sqrt 5)/2,-(3 (470 \pm 213 \sqrt 5))^3(1 \pm \sqrt 5)/2\\
 &(-1/2,3/8)&\text{non}      &-&\Q(\sqrt{-15},\sqrt{177},\sqrt{1257})\\
& (-1/2,-3/8) &\text{yes}      &-12,-48&2\cdot 30^3,4 (15 (30 \pm 17 \sqrt 3))^3\\
 &(-4/3,143/27) &\text{non}      &-&\Q(\sqrt{37},\sqrt{-143},\sqrt{2077})\\
&(-4/3,-143/27) &\text{yes}& -332&\Q(\sqrt 3,\sqrt{31})\\  \hline\hline
 230 &\infty'& \text{yes}   &  -40   &(6 (65 \pm 27 \sqrt 5))^3\\
    &(0,  1) &\text{yes}&-20&(2 (25 \pm  13 \sqrt 5))^3 \\
  &(0,  -1) &\text{yes}&-15&-(3 (25 \pm 9 \sqrt 5)/2 )^3 (-1 \pm \sqrt 5)/2\\
    &(1,  5) &\text{yes}&-520 &\Q(\sqrt 5,\sqrt{13})\\
  &(1,  -5) &\text{yes}&-120&\Q(\sqrt 2, \sqrt 5 )\\
  &(-2,  5) &\text{yes}&-15, -60&-(3 (25 \pm 9 \sqrt 5)/2 )^3 (-1 \pm \sqrt 5)/2,-(3 (470 \pm 213 \sqrt 5))^3(1 \pm \sqrt 5)/2\\
  &(-2,  -5) &\text{yes}&-1380 &  \Q(\sqrt 3,\sqrt 5, \sqrt {23})\\
     &(3,35) &\text{non }&-&\Q(\sqrt{685},\sqrt{705},\sqrt{19043})\\
      &(3,-35) &\text{yes}&-180&\Q(\sqrt{3},\sqrt 5 )\\
      \hline\hline
266 &\infty'&\text{yes}&-52& (30 (31 \pm  9 \sqrt{13}))^3 \\
&(-1, 1) &\text{yes} &-84& \Q(\sqrt 3, \sqrt 7)\\
&(-1, -1) &\text{yes}&-3,-12 &0, 2\cdot 30^3 \\
&(0, 1) &\text{yes}&-280& \Q(\sqrt 2, \sqrt 5)\\
&(0, -1) &\text{yes}&-40&(6 (65 \pm 27 \sqrt 5))^3 \\
 &(-5/2, 83/8) &\text{non}&-& \Q(\sqrt{1041},\sqrt{-415},\sqrt{105})\\
 & (-5/2, -83/8) &\text{yes}&-532& \Q(\sqrt 17,\sqrt {19})\\
  \hline\hline
  285 &\infty'&\text{yes}&-51&-(48 (37 \pm 9 \sqrt{17}))^3 (-4 \pm \sqrt{17})  \\
&(-1, 4) &\text{yes} &-15& -(3 (25 \pm 9 \sqrt 5)/2 )^3 (-1 \pm \sqrt 5)/2\\
&(-1, -4) &\text{yes}&-60& -(3 (470 \pm 213 \sqrt 5))^3(1 \pm \sqrt 5)/2\\
&(0, 0) &\text{yes}&-3,-75& 0, -(48 (-69 \pm  31 \sqrt 5))^3 (\pm \sqrt 5)\\
&(3, 24) &\text{non}&-&\Q (\sqrt{3},\sqrt{95},\sqrt{60197})\\
&(3, -24) &\text{yes}&-240& \Q(\sqrt 3 ,\sqrt 5 )\\
 &(-3/2, 57/8) &\text{non}&& \Q(\sqrt{-79},\sqrt{57},\sqrt{11985})\\
 & (-3/2, -57/8) &\text{yes}&-1995&\Q( \sqrt 5,\sqrt{21},\sqrt{57}) \\
  \hline\hline
   286 &\infty'&\text{yes}&-40&(6 (65 \pm 27 \sqrt 5))^3 \\
&(-1, 4) &\text{yes} &-52&(30 (31 \pm  9 \sqrt{13}))^3  \\
&(-1, -4) &\text{yes}&-88& (60 (155 \pm  108 \sqrt 2))^3 \\
&(5/2, 143/8) &\text{non}&-&\Q (\sqrt{39},\sqrt{168917},\sqrt{232})\\
&(5/2, -143/8) &\text{non}&-& \Q(\sqrt{1841},\sqrt{-3367},\sqrt{37609})\\ \hline\hline
 357 &\infty'&\text{yes}&-168& \Q(\sqrt 6 ,\sqrt{14})\\
&(-1, 4) &\text{non} &-& \Q(\sqrt{293},\sqrt{89997},\sqrt{21})\\
&(-1, -4) &\text{yes}&-35& -(16(15 \pm  7 \sqrt 5))^3 \\ \hline\hline
\end{array}
$$
\end{small}

\subsection{ $N$ is a product of four primes: $\Gal (\Q( j)/\Q)\hookrightarrow(\Z/2\Z)^4$.}
\begin{small}
$$\begin{array}{c|c|c|c|c|}
N&  \text{point} & \text{CM}  & D&  j\text{ or }\Q(j)\\ \hline\hline
   390&\infty'&\text{yes}&-5460& \Q(\sqrt 3 , \sqrt 5 ,\sqrt 7,\sqrt{11})\\
&(0, 1) &\text{yes} &-120&\Q(\sqrt 2 ,\sqrt 5) \\
&(0, -1) &\text{yes}&-420& \Q(\sqrt 3 , \sqrt 5,\sqrt 7  )\\
&(1, 2) &\text{yes} &-660&\Q(\sqrt 3, \sqrt  5 ,\sqrt {11} ) \\
&(1,-2) &\text{yes}&-4,-36&12^3, -(4 (102 \pm  61 \sqrt 3))^3 (-2 \pm  \sqrt 3)\\ \hline\hline
\end{array}
$$
\end{small}


\vskip 0.4 cm

\noindent{Francesc Bars Cortina}\\
{Departament Matem\`atiques, Edif. C, Universitat Aut\`onoma de Barcelona\\
08193 Bellaterra, Catalonia}\\
{francesc@mat.uab.cat}

 \vspace{1cm}

\noindent
{Josep Gonz\'alez Rovira}\\
{Departament de Matem\`atiques, Universitat Polit\`ecnica de Catalunya EPSEVG,\\
Avinguda V\'{\i}ctor Balaguer 1, 08800 Vilanova i la Geltr\'u,
Catalonia}\\
{josep.gonzalez@upc.edu}

\vspace{1cm}

\noindent{Xavier Xarles}\\
{Departament Matem\`atiques, Edif. C, Universitat Aut\`onoma de Barcelona\\
08193 Bellaterra, Catalonia}\\
{xarles@mat.uab.cat}
\end{document}